\documentclass[reqno]{amsart}\usepackage{amsmath}

\usepackage{amsfonts}
\usepackage{amssymb}
\usepackage{hyperref}

\begin{document}
\title[Periodic Boundary
Value Problems]{Remarks on the periodic boundary
value problems for nonlinear second order ordinary differential
equations}
\author[F. Haddouchi, S. Benaicha]{Faouzi Haddouchi, Slimane Benaicha}
\address{Faouzi Haddouchi\\
Department of Physics, University of Sciences and Technology
of Oran-MB, El Mnaouar, BP 1505, 31000 Oran, Algeria\\
Laboratory of Fundamental and Applied Mathematics of Oran (Lmfao), Department of Mathematics, University of Oran 1 Ahmed Benbella, 31000 Oran,
Algeria}
\email{fhaddouchi@gmail.com}
\address{Slimane Benaicha \\
Laboratory of Fundamental and Applied Mathematics of Oran (Lmfao), Department of Mathematics, University of Oran 1 Ahmed Benbella, 31000 Oran,
Algeria} \email{slimanebenaicha@yahoo.fr}

\begin{abstract}
This paper is devoted to study the existence of solutions and the
monotone method of second-order periodic boundary value problems
when the lower and
upper solutions $\alpha$ and $\beta$ violate the boundary conditions $%
\alpha(0)=\alpha(2\pi)$ and $\beta(0)=\beta(2\pi)$. We present
several comparison results. We show that the method of lower and
upper solutions coupled with the monotone iterative technique is
valid to obtain constructive proof of existence of solutions.
\end{abstract}

\maketitle \numberwithin{equation}{section}
\newtheorem{theorem}{Theorem}[section]
\newtheorem{lemma}[theorem]{Lemma} \newtheorem{proposition}[theorem]{%
Proposition} \newtheorem{corollary}[theorem]{Corollary} \newtheorem{remark}[%
theorem]{Remark}
\newtheorem{exmp}{Example}[section]

\section{Introduction}
It is well-known that the upper and lower solution method and the iterative technique
is a powerful tool for proving the existence results for nonlinear boundary value problem(BVP for short), see \cite{MaYang,Cab,Nieto}.
In \cite{Nieto2Alva}, the authors introduced a new concept of upper and lower solution and then they
proved the existence results for nonlinear first order boundary value problem.
In this paper, we are concerned with the existence of solutions and
the monotone method for periodic boundary value problems (PBVP, for
short) of second-order differential equations of the form
\begin{equation} \label{eqn1}
\text{\ }-{u^{\prime \prime }}=f(t,u),\text{ \ }u(0)=u(2\pi ),\text{ }{%
u^{\prime }}(0)={u^{\prime }}(2\pi ),
\end{equation}
where $f\in C[I\times \mathbb{R},\mathbb{R}]$, $I=[0,2\pi ]$.
It is well known that if $\alpha $ and $\beta $ are lower and upper
solutions for \eqref{eqn1}, respectively, with $\alpha \leq \beta $ on $%
[0,2\pi ]$ and satisfy%
\begin{equation}  \label{eqn2}
\alpha (0)=\alpha (2\pi )\text{ , \ }\beta (0)=\beta (2\pi )
\end{equation}%
\begin{equation}  \label{eqn3}
{\alpha ^{\prime }}(0)\geq {\alpha ^{\prime }}(2\pi )\text{ , \
}{\beta ^{\prime }}(0)\leq {\beta ^{\prime }}(2\pi ),
\end{equation}%
then the existence of solutions has been proved by using an abstract
existence theorem (see, for instance \cite{Kan1Lakshm}, \cite{Kan2Lakshm}).
Moreover, the extreme solutions have been obtained as limits of
monotone iterates.
The problem of proving the existence of solutions when the conditions %
\eqref{eqn3} are violated has been investigated in detail; (see,
e.g.,\cite{Lakshm1Lee}, \cite{Lakshm2} for an excellent
bibliography). In \cite{Lakshm1Lee}, \cite{Lakshm2}, only the
monotone iterative method has been discussed. In this paper, our
main objective is to investigate the case when $\alpha $ and $\beta
$ violate the relations of \eqref{eqn2}. Motivated by the studies in
\cite{Lakshm1Lee}, \cite{Lakshm2}, and \cite{Nieto2Alva}, we will establish some
comparison results that will be useful later. Also, we will show the validity of the upper and lower solutions method as well as the monotone iterative technique if adding one-side Lipschitz condition to $f(t,u)$. Finally, we give an illustrative example showing the usefulness of our main results.

\section{Linear periodic boundary value problem}

We need the following lemma for the sequel. The proof of this lemma
is omitted as it can be obtained by direct computations.

\begin{lemma}\label{Lemma 2.1}
Let $M$, $\lambda$, $\mu \in \mathbb{R}$ with $M\neq0$ and $\sigma
\in C([0,2\pi ],\mathbb{R})$. The solution of the linear problem
\begin{equation}  \label{eqn4}
\left\{
\begin{array}{l}
-{u^{\prime\prime}}(t)+M^{2}u(t)=\sigma (t) \\
u(0)-u(2\pi )=\mu \\
{u^{\prime}}(0)-{u^{\prime}}(2\pi )=\lambda%
\end{array}%
\right.
\end{equation}
is given by
\begin{equation}  \label{eqn5}
u(t)=C_{1}e^{Mt}+C_{2}e^{-Mt}-\frac{e^{Mt}}{2M}\int_{0}^{t}e^{-Ms}\sigma
(s)ds+\frac{e^{-Mt}}{2M}\int_{0}^{t}e^{Ms}\sigma (s)ds
\end{equation}
where
\[
C_{1}=\frac{1}{2(1-e^{2\pi M})}(\mu +\frac{\lambda }{M})-\frac{e^{2\pi M}}{%
2M(1-e^{2\pi M})}\int_{0}^{2\pi }e^{-Ms}\sigma (s)ds,
\]
\[
C_{2}=\frac{1}{2(1-e^{-2\pi M})}(\mu -\frac{\lambda }{M})+\frac{e^{-2\pi M}}{%
2M(1-e^{-2\pi M})}\int_{0}^{2\pi }e^{Ms}\sigma (s)ds.
\]
\end{lemma}

\begin{theorem}\label{Theorem 2.2}
Let $u\in C^{2}([0,2\pi ], \mathbb{R})$ and suppose that there
exists $ \omega\in C([0,2\pi ],\mathbb{R})$, $\omega \geq 0$ on
$I=[0,2\pi ]$ such that
\begin{equation} \label{eqn6}
\left\{
\begin{array}{l}
-{u^{\prime\prime}}(t)+M^{2}u(t)+ \omega(t)\leq 0\text{, }t\in I \\
u(0)-u(2\pi )=\mu \\
{u^{\prime}}(0)-{u^{\prime}}(2\pi )=\lambda%
\end{array}%
\right.
\end{equation}
and%
\begin{equation} \label{eqn7}
\int_{0}^{2\pi }(e^{Ms}+e^{-M(s-2\pi )})\omega(s)ds\geq \mu
M(e^{2\pi M}-1)-\lambda (e^{2\pi M}+1).
\end{equation}
Then $u\leq 0$ on $I$.
\end{theorem}

\begin{proof}
Note that \ $-{u^{\prime \prime }}(t)+M^{2}u(t)+\omega(t)=\sigma $
with $\sigma
\leq 0$ on $I$. Using \eqref{eqn5} and \eqref{eqn7} we have that%
\begin{align*}
u(0) &=C_{1}+C_{2} \\
&=\frac{1}{2(1-e^{2\pi M})}(\mu +\frac{\lambda
}{M})+\frac{1}{2(1-e^{-2\pi
M})}(\mu -\frac{\lambda }{M})- \\
&\int_{0}^{2\pi }\frac{e^{Ms}+e^{-M(s-2\pi )}}{2M(1-e^{2\pi
M})}\sigma (s)ds+\int_{0}^{2\pi }\frac{e^{Ms}+e^{-M(s-2\pi
)}}{2M(1-e^{2\pi M})}\omega(s)ds
\end{align*}
and
\begin{align*}
{\small u(0)}&\le \frac{1}{{\small 2M(1-e}^{2\pi M}{\small
)}}\left[ \int_{0}^{2\pi }{\small (e}^{Ms}{\small +e}^{-M(s-2\pi
)}{\small )\omega(s)ds-\mu
M(e}^{2\pi M}{\small -1)+\lambda (e}^{2\pi M}{\small +1)}\right] \\
&\le 0.
\end{align*}%
\
Now, if $u(t)>0$ for some $t\in \left( 0,2\pi \right] $, then $\max \left\{ u(t),\text{ }t\in \left[ 0,2\pi \right] \right\} >0$. By taking $u(t^{\star })=\max \{ u(t),\text{ }t\in [ 0,2\pi]\}$, we have ${u^{\prime}}(t^{\star })=0$ \ and ${u^{\prime\prime}}(t^{\star })\leq 0$.

From the inequality $\ -{u^{\prime\prime}}(t^{\star
})+M^{2}u(t^{\star })+\omega(t^{\star })\leq 0$ \ we get that
${u^{\prime\prime}}(t^{\star })>0$ which is a contradiction.
\end{proof}

This result gives us several useful consequences.

\begin{corollary}\label{Corollary 2.3}
Let $u\in C^{2}([0,2\pi ],\mathbb{R})$ such that $-{u^{\prime\prime}}(t)+M^{2}u(t)\leq 0$ for $t\in I$ and%
\[
u(0)-u(2\pi )\leq \frac{1}{M\tanh (\pi
M)}({u^{\prime}}(0)-{u^{\prime}}(2\pi )).
\]
Then $u\leq 0$ on $I.$
\end{corollary}
\begin{proof}
It is result immediately from Theorem \ref{Theorem 2.2}.
\end{proof}

\begin{corollary}\label{Corollary 2.4}
Let $u\in C^{2}([0,2\pi ],\mathbb{R})$, $\omega \in \mathbb{R}$,
$\omega \geq 0$ such that $-{u^{\prime\prime}}(t)+M^{2}u(t)\leq
-\omega $ \ for $t\in
I $,%
\[
u(0)-u(2\pi )>\frac{1}{M\tanh (\pi M)}({u^{\prime}}(0)-{u^{\prime}}(2\pi ))%
\text{ }
\]
and
\begin{equation}\label{eqn8}
\omega \geq \frac{M^{2}}{2}[u(0)-u(2\pi )-\frac{1}{M\tanh (\pi M)}({u^{\prime}}(0)-%
{u^{\prime}}(2\pi ))].
\end{equation}%
Then $u\leq 0$ on $I$.
\end{corollary}

\begin{proof}
We have \
\[
\int_{0}^{2\pi }(e^{Ms}+e^{-M(s-2\pi )})\omega
ds=\frac{2\omega}{M}(e^{2\pi M}-1).
\]
From \eqref{eqn8}, we obtain
\begin{align*}
\int_{0}^{2\pi }(e^{Ms}+e^{-M(s-2\pi )}) \omega ds &\geq \frac{2}{M}(e^{2\pi M}-1)%
\frac{M^{2}}{2}[\mu -\frac{1}{M\tanh (\pi M)}\lambda ] \\
&\geq \mu M(e^{2\pi M}-1)-\lambda (e^{2\pi M}+1)\text{.}
\end{align*}
By virtue of Theorem \ref{Theorem 2.2} we get that $u\leq 0$ on $I$.
\end{proof}

\begin{theorem}\label{Theorem 2.5}
Let $u\in C^{2}([0,2\pi ],\mathbb{R})$ and suppose that there exists $\omega \in C([0,2\pi ],\mathbb{R})$, $\omega \leq 0$ on $I=[0,2\pi ]$ such that%
\begin{equation}\label{eqn9}
\left\{
\begin{array}{l}
-{u^{\prime\prime}}(t)+M^{2}u(t)+\omega(t)\geq 0\text{, }t\in I \\
u(0)-u(2\pi )=\mu \\
{u^{\prime}}(0)-{u^{\prime}}(2\pi )=\lambda%
\end{array}%
\right.
\end{equation}
and%
\begin{equation}\label{eqn10}
\int_{0}^{2\pi }(e^{Ms}+e^{-M(s-2\pi )})\omega (s)ds\leq \mu
M(e^{2\pi M}-1)-\lambda (e^{2\pi M}+1).
\end{equation}
Then $u\geq 0$ on $I$.
\end{theorem}

\begin{proof}
The same argument as in Theorem \ref{Theorem 2.2} will be used.
\end{proof}
\section{Existence theorems}

We now consider the nonlinear boundary value problem%
\begin{equation}  \label{eqn11}
-{u^{\prime \prime }}(t)=f(t,u(t)),\text{ \ }u(0)=u(2\pi ),\text{ }{%
u^{\prime }}(0)={u^{\prime }}(2\pi )
\end{equation}
where $f\in C([0,2\pi ]\times \mathbb{R},\mathbb{R})$. Relative to
the lower and upper solutions $\alpha$, $\beta$ of the problem
\eqref{eqn11}, we shall list the following assumptions:
\begin{itemize}

\item[(F1)]$\alpha, \beta \in C^{2}([0,2\pi ],\mathbb{R})$, $%
\alpha(t)\leq\beta(t)$ on $\left[0,2\pi\right]$, and there exists $M\neq0$ such that

$f(t,u)-f(t,v)\geq-M^{2}(u-v)$, $t\in\left[0,2\pi\right]$ whenever $\alpha(t)\leq v\leq u\leq
\beta(t)$.

\item[(F2)] $-{\alpha^{\prime \prime }}\leq f(t,\alpha)$, $t\in\left[0,2\pi%
\right]$, $\alpha(0)-\alpha(2\pi)<0$, ${\alpha^{\prime}}(0)-{\alpha^{\prime}}%
(2\pi)\geq0$ and \newline
$-{\beta^{\prime \prime }}\geq f(t,\beta)$, $t\in\left[0,2\pi\right]$, $%
\beta(0)-\beta(2\pi)>0$, ${\beta^{\prime}}(0)-{\beta^{\prime}}(2\pi)\leq0$.%
\end{itemize}

\begin{theorem}\label{Theorem 3.1}(Upper and lower solutions method)
Let {\rm (F1)} and {\rm (F2)} hold. Then, the problem \eqref{eqn11} is
solvable.
\end{theorem}
\begin{proof}
For $u\in \mathbb{R}$, we set $p(t,u)=\min \{\beta (t),\max
\{u,\alpha (t)\}\}$. Thus, for $u\in C^{2}([0,2\pi ],\mathbb{R} )$ \
we define the function $\hat{u}(t)=p(t,u(t))$. We now consider the
modified problem

\begin{equation}  \label{eqn12}
\left\{
\begin{array}{l}
-{u^{\prime\prime}}(t)+M^{2}u(t)=F(t,u(t))\text{, }t\in I \\
u(0)=u(2\pi ),\text{ }{u^{\prime}}(0)={u^{\prime}}(2\pi )%
\end{array}%
\right.
\end{equation}
where \[F(t,u(t))=f(t,p(t,u(t)))+M^{2}p(t,u(t)).\]
The problem \eqref{eqn12} admits a solution given by \eqref{eqn5} where $\mu =\lambda =0,$ and $\sigma (t)=F(t,u(t))$. Note that if $u$ is solution of %
\eqref{eqn12} such that $\alpha \leq u\leq \beta $ on $I$, then $u$
is actually a solution of \eqref{eqn11}. We shall prove that any solution of %
\eqref{eqn12} is such that $\alpha \leq u\leq \beta $ on $I$. Thus,
we obtain that \eqref{eqn11} has at least one solution.

Indeed, let $v=\alpha -u$. Then using  {\rm (F1)} we have that for
every $t\in I$,
\begin{align*}
-{v^{\prime\prime}}+M^{2}v+({\alpha^{\prime\prime}}+f(t,\alpha ))&=-[f(t,%
\hat{u})-f(t,\alpha )+M^{2}(\hat{u}-\alpha )] \\
&\le 0.
\end{align*}
Now, using Theorem \ref{Theorem 2.2}, we obtain that $v\leq 0$ on $I.$
Similarly, we have that $u\leq \beta $ on $I$.
\end{proof}

\begin{theorem}\label{Theorem 3.2}(Monotone method).
Let  {\rm (F1)} and  {\rm (F2)} hold. Then, there exist monotone
sequences $\left\{ \alpha _{n}\right\} \nearrow \phi $, and $\left\{
\beta _{n}\right\} \searrow \psi $ uniformly on $I$ with $\alpha
_{0}=\alpha $ and $\beta _{0}=\beta $. Here $\phi $ and $\psi $ are
the minimal and maximal solutions of \eqref{eqn11} respectively
between $\alpha $ and $\beta $, that is, if $u$ is a solution of
\eqref{eqn11} with $\alpha \leq u\leq \beta $ on $I$, then $\phi
\leq u\leq \psi $ on $I$. Moreover, these sequences are such that
$\alpha _{0}\leq ...\leq \alpha _{n}\leq ...\leq \beta _{m}\leq
...\leq \beta _{0}$, for every $n$, $m\in \mathbb{N}$.
\end{theorem}

\begin{proof}
For $\eta \in \left[ \alpha ,\beta \right]=\left\{ \eta \in C([0,2\pi ],%
\mathbb{R}):\alpha \leq \eta \leq \beta \text{ on }I\right\} $, let
us consider the linear periodic boundary value problem%
\begin{equation}  \label{eqn13}
\left\{
\begin{array}{l}
-{u^{\prime\prime}}(t)+M^{2}u(t)=\sigma _{\eta }(t) \\
u(0)=u(2\pi ),\text{ }{u^{\prime}}(0)={u^{\prime}}(2\pi )%
\end{array}%
\right.
\end{equation}

where \[\sigma (t)=\sigma _{\eta }(t)= f(t,\eta (t))+M^{2}\eta
(t).\]
This problem has a unique solution given by%
\begin{equation}  \label{eqn14}
u(t)=C_{1}e^{Mt}+C_{2}e^{-Mt}-\frac{e^{Mt}}{2M}\int_{0}^{t}e^{-Ms}\sigma
(s)ds+\frac{e^{-Mt}}{2M}\int_{0}^{t}e^{Ms}\sigma (s)ds,
\end{equation}

where%
\begin{eqnarray*}
C_{1} &=&\text{$-\frac{e^{2\pi M}}{2M(1-e^{2\pi M})}\int_{0}^{2\pi
}e^{-Ms}\sigma (s)ds,$} \\
\text{ }C_{2} &=&\frac{e^{-2\pi M}}{2M(1-e^{-2\pi M})}\int_{0}^{2\pi
}e^{Ms}\sigma (s)ds.
\end{eqnarray*}
To prove uniqueness, let $v(t)$ another solution of \eqref{eqn13}
and define
$w(t)=u(t)-v(t)$. We see that%
\[
-{w^{\prime\prime}}(t)+M^{2}w(t)=0,\text{ }w(0)=w(2\pi ),\text{ }{w^{\prime}}%
(0)={w^{\prime}}(2\pi ).
\]

Hence, by Corollary \ref{Corollary 2.3}, it follows that $w(t)\equiv 0$,
which shows $v(t)=$ $u(t)$. Hence for any $\eta \in \lbrack \alpha
,\beta ]$, we define an operator $A$ by $A\eta =u$, where $u$ is the unique solution of %
\eqref{eqn13}. The operator $A$ is well defined from $[\alpha ,\beta ]$ to $%
[\alpha ,\beta ]$ and $A$ is monotone nondecreasing on $[\alpha
,\beta ]$. Indeed, let $\eta \in \lbrack \alpha ,\beta ]$ and define
$v=\alpha -u$.
Thus, for all $t\in I$ we have%
\begin{align*}
-{v^{\prime\prime}}(t)+M^{2}v(t)+({\alpha^{\prime\prime}}+f(t,\alpha ))&=-%
\left[ f(t,\eta )-f(t,\alpha )+M^{2}(\eta -\alpha )\right] \\
&\le 0\text{.}
\end{align*}
Hence, by Theorem \ref{Theorem 2.2}, it follows that $v\leq 0$, which shows $%
\alpha \leq u$ on $I$. Similarly, we get that $u\leq \beta $ on $I$.

To show the monotonicity of $A$, let $\eta _{1}$, $\eta _{2}\in
\left[ \alpha ,\beta \right]$ such that $\eta _{1}\leq \eta _{2}$.
Let $A\eta_{1}=u_{1}$ and $A\eta _{2}=u_{2}$. Setting $w=u_{1}-u_{2}$ and using  {\rm (F1)} , we get%
\[
-{w^{\prime\prime}}(t)+M^{2}w(t)\leq 0\text{, }w(0)=w(2\pi ),\text{ }{%
w^{\prime}}(0)={w^{\prime}}(2\pi ).
\]
Hence, by Corollary \ref{Corollary 2.3}, we have $w(t)\leq 0$ and this
implies that $A$ is monotone on $[\alpha ,\beta ]$.

We now define $\{\alpha _{n}\}$,  $\{\beta _{n}\}$ with $\alpha _{0}=\alpha $%
,  $\beta _{0}=\beta $ by%
\[
\alpha _{n+1}=A\alpha _{n}\text{, \ }\beta _{n+1}=A\beta _{n}.\]
Then we have%
\[
\alpha =\alpha _{0}\leq ...\leq \alpha _{n}\leq ...\leq \beta
_{n}\leq ...\leq \beta _{0}=\beta ,\text{ \ \ \ \ \ }t\in I.
\]
It then follows, by using standard argument \cite{Lakshm1Lee}, that $\underset{%
n\longrightarrow \infty }{\lim }\alpha _{n}(t)=\phi (t)$ and $\underset{%
n\longrightarrow \infty }{\lim }\beta _{n}(t)=\psi (t)$ uniformly
and monotonically on $I$, and that $\phi $ and $\psi $ are the
minimal and maximal solutions of \eqref{eqn11}, respectively. This
completes the proof.
\end{proof}

\section{Application}
Now, we give two examples to illustrate our results.
\begin{exmp}\label{exmp 4.1}
Consider the following BVP
\begin{equation}\label{eqn15}
-{u^{\prime \prime }}(t)+2u(t)=t, \  \ u(0)=u(2\pi), \  \ {u^{\prime  }}(0)={u^{\prime }}(2\pi).
\end{equation}
We define the function $f:[0,2\pi]\times \mathbb{R}\longrightarrow\mathbb{R}$ by
\[f(t,u)=t-2u, \  \ t\in[0, 2\pi]. \]
In this case, for every $u, v \in\mathbb{R}$ with $v\leq u$ and $t\in[0, 2\pi]$, we have
\[f(t,u)-f(t,v)=-2(u-v)\geq-M^{2}(u-v), \text{where}\ M=2.\]
Now, let $\alpha(t)=t-2\pi$ and $\beta(t)=2\pi-\frac{1}{4}t$. Thus,
\[\alpha(0)-\alpha(2\pi)=-2\pi<0, \ \  \alpha^{\prime}(0)-\alpha^{\prime}(2\pi)=0, \]
and for every $t\in[0, 2\pi]$ we have
\[-\alpha^{\prime \prime }\equiv0\leq f(t,\alpha(t))=4\pi-t. \]
On the other hand, by the definition of $\beta(t)$, we have
\[\beta(0)-\beta(2\pi)=\frac{\pi}{2}>0, \ \  \beta^{\prime}(0)-\beta^{\prime}(2\pi)=0, \]
and for every $t\in[0, 2\pi]$ we have
\[-\beta^{\prime \prime }\equiv0\geq f(t,\beta(t))=\frac{3}{2}t-4\pi. \]
Finally, by the choice of $\alpha(t)$ and $\beta(t)$ we see that $\alpha(t)\leq \beta(t)$.
So, according to Theorem \ref{Theorem 3.1}, we can conclude that there exists a solution $u$ of \eqref{eqn15} with $\alpha \leq u \leq \beta$.
\end{exmp}

\begin{exmp}\label{exmp 4.2}
Consider the following BVP
\begin{equation}\label{eqn16}
-{u^{\prime \prime }}(t)+\frac{1}{4\pi^{2}}u^{2}(t)=t+\pi, \  \ u(0)=u(2\pi), \  \ {u^{\prime  }}(0)={u^{\prime }}(2\pi),
\end{equation}
where $f(t,u)=-\frac{1}{4\pi^{2}}u^{2}(t)+t+\pi$. 

Now, let $\alpha(t)=t-2\pi$. Thus,
\[\alpha(0)-\alpha(2\pi)=-2\pi<0, \ \  \alpha^{\prime}(0)-\alpha^{\prime}(2\pi)=0\geq0, \]
and for every $t\in[0, 2\pi]$, we have by a simple calculation
\[-\alpha^{\prime \prime }\equiv0\leq f(t,\alpha(t))=-\frac{1}{4\pi^{2}}t^{2}+(\frac{1}{\pi}+1)t+\pi-1. \]
Analogously, let $\beta(t)=9\pi-t$. Thus,
\[\beta(0)-\beta(2\pi)=2\pi>0, \ \  \beta^{\prime}(0)-\beta^{\prime}(2\pi)=0\leq0,\]
and for every $t\in[0, 2\pi]$ we have
\[-\beta^{\prime \prime }\equiv0\geq f(t,\beta(t))=-\frac{1}{4\pi^{2}}t^{2}+(1+\frac{9}{2\pi})t+\pi-\frac{81}{4}. \]
By the choice of $\alpha(t)$ and $\beta(t)$ we see that $\alpha(t)\leq \beta(t)$ for all $t\in[0, 2\pi]$.

Finally, for every $t\in[0, 2\pi]$ with $\alpha(t)\leq v\leq u\leq v\leq \beta(t)$, we have
\begin{eqnarray*}
f(t,u)-f(t,v)&=&-\frac{1}{4\pi^{2}}(u+v)(u-v)\\
&\geq& \frac{1}{4\pi^{2}}(-2\beta)(u-v)\\
&\geq& \frac{1}{4\pi^{2}}(-18\pi)(u-v)\\
&=&-\Big(\frac{3}{\sqrt{2\pi}}\Big)^{2}(u-v).
\end{eqnarray*}
All conditions of Theorem \ref{Theorem 3.1} are satisfied. Thus, there exists a solution $u$ of \eqref{eqn16} with $\alpha \leq u \leq \beta$.
\end{exmp}

\end{document}